\documentclass[]{article}
\usepackage[english]{babel}
\usepackage[normalem]{ulem}
\usepackage[all]{xy}
\usepackage{amsmath,amsthm,booktabs,color,dsfont,epsfig,graphicx,mathrsfs,multirow,scalefnt}
\usepackage{amsfonts, amssymb}
\usepackage{float}
\usepackage{caption}
\usepackage{xcolor}

\usepackage{mathtools}

\usepackage{xr}

\usepackage[most]{tcolorbox}

\tcbset{
  highlightstyle/.style={
    colback=gray!7,     
    colframe=gray!7,    
    boxrule=0.4pt,       
    arc=2pt,             
    left=4pt, right=4pt, top=2pt, bottom=2pt, 
    enhanced             
  }
}

\newtheorem{theo}{Theorem}[section]
\newtheorem*{theo*}{Theorem}
\newtheorem{cor}[theo]{Corollary}
\newtheorem{lem}[theo]{Lemma}
\newtheorem{claim}[theo]{Claim}
\newtheorem{prop}[theo]{Proposition}

\newcommand{\Z}{\mathbb{Z}}
\newcommand{\N}{\mathbb{N}}
\newcommand{\NZ}{\mathbb{M}}

\newcommand{\s}{\sigma}

\newcommand{\x}{\mathbf{x}}
\newcommand{\y}{\mathbf{y}}

\newcommand{\bfN}{\mathbf{N}}
\newcommand{\bfM}{\mathbf{M}}

\title{Corrigendum to the paper ‘Some notes on the classification of shift spaces: Shifts of Finite Type; Sofic Shifts; and Finitely Defined Shifts’ [Bulletin of the Brazilian Mathematical Society, New Series (2022), 53, 981-1031].}

\author{
\small{Marcelo Sobottka}\\
\footnotesize{UFSC -- Department of Mathematics}\\
\footnotesize{88040-900 Florian\'{o}polis - SC, Brazil}\\
\footnotesize{\texttt{marcelo.sobottka@ufsc.br}}
}

\date{}

\begin{document}

\maketitle

\begin{abstract} This paper is a corrigendum to the article ‘Some notes on the classification of shift spaces: Shifts of Finite Type; Sofic Shifts; and Finitely Defined Shifts’. In this article we correct Lemma 5.3. Therefore, we follow correcting statements and proofs of subsequent results that depend on Lemma 5.3. 
\end{abstract}


\vspace{0.5cm}
\tiny

\begin{center}

\begin{minipage}[t]{.9\textwidth}
\hrule
\colorbox{gray!10}{%
    \begin{minipage}{.985\textwidth}

\noindent This preprint has not undergone peer review or any post-submission improvements or corrections. The Version of Record of this article is published in Bulletin of the Brazilian Mathematical Society, New Series, and shall be cited as:\\

{\bf Sobottka, M.  Correction: Some notes on classification of shift spaces: Shifts of Finite Type; Sofic Shifts; and Finitely Defined Shifts. Bull Braz Math Soc, New Series 56, 50 (2025). https://doi.org/10.1007/s00574-025-00470-7.}

\end{minipage}%
}
\hrule
\end{minipage}

\end{center}
\normalsize

\section{Introduction}\label{sec:Introduction}

The aim of this paper is to present a correct version for Lemma \ref{lem:SFTs_finite-tersect} in \cite{Sobottka2022}. Such a lemma in its original (and incorrect) version stated that shifts of finite type on the lattices $\N$ or $\Z^d$ (and independently of the cardinality of the alphabet) always present the property of nested cylinders, that is, they would be such that if $\{W_\ell\}_{\ell\geq 1}$ is a family of nested non-empty cylinders, then $\bigcap_{\ell= 1}^\infty W_\ell\neq\emptyset$. In fact, whenever the alphabet is finite and the lattice is countable, the compactness implies that such property holds independently on the type of the shift. However, if the alphabet is infinite, it is possible that a shift space has not such a property even if the lattice is $\Z$. Finding sufficient conditions for the property of nested cylinders is extremely important for characterizing continuous shift-commuting maps between shift spaces (see Theorem \ref{theo: image_of_shifts-intersection} and Proposition \ref{prop:higher_block_shift-NZ_countable} in \cite{Sobottka2022}).\\

The idea behind the original proof of Lemma \ref{lem:SFTs_finite-tersect} was to fix the entries gradually, relying on the fact that the forbidden words had finite length. This would yield an increasing sequence of patterns that should lead to the existence of a point in the intersection of all cylinders. Unfortunately, such an argument has a flaw: the way the nested cylinders fixed patterns left gaps that, in each cylinder, had to be filled in a different way. 

However, for the specific case of Markov shifts on the lattices $\N$ or $\Z$, the argument used in that proof works. So we will present in this article a correct statement and proof for Lemma \ref{lem:SFTs_finite-tersect}. It is worth to note that Lemma \ref{lem:SFTs_finite-tersect} was used in several results in \cite{Sobottka2022}, and hence, a more restrictive version of it will imply that those results should also be corrected. Figure 1.E below presents the chain of results based on Lemma \ref{lem:SFTs_finite-tersect} in the original paper.

\begin{equation*}\hspace{-2cm}
\xymatrix{
  &       *+[F]{Lemma\ 5.3 }\ar[r] \ar[d]       &   *+[F]{Corollary\ 5.4}      &      \\
*+ [F]{Corollary\ 5.10}\ar[d]   &    *+[F]{Proposition\ 5.8.(C2)}\ar[r]\ar[l]         &  *+[F]{Proposition\ 6.5 }\ar[d]^{using\ Prop.\ 5.8.(C2)}         &    \\
{ \begin{array}{c}Section\ 7\\
   (except\ Theorem\ 7.8)\end{array}}   &        &  *+[F]{ Corollary\ 6.6 }\ar[lld]\ar[d]\ar[rd]    &  \\
 Theorems\ 7.5,\ 7.6,\ 7.10  &                 &  *+[F]{ Theorem\ 7.4 }\ar[d]\ar[r]        &   *+[F]{ Theorem\ 7.9} \\
   &            &  *+[F]{ Claim\ 9.4 }         &  
  }
  \end{equation*}
  \captionsetup{labelformat=empty} 

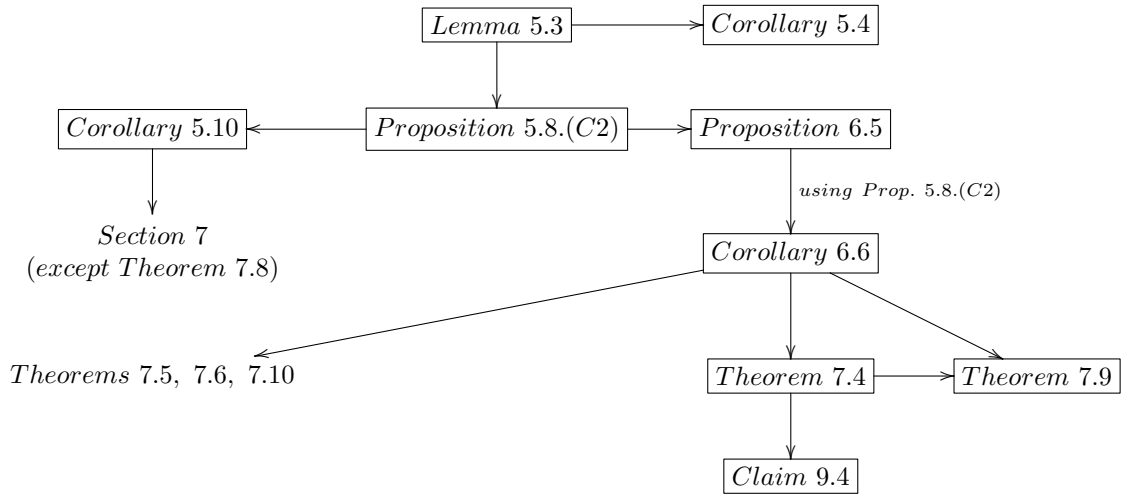
\captionof{figure}{\textbf{Figure 1.E:} The scheme above shows how the results are related to Lemma 5.3 in the original article. A box around a result indicates a strong dependence on that lemma.}

\vspace{0.75cm}

The change in Lemma \ref{lem:SFTs_finite-tersect} caused significant changes in the generality of Corollary \ref{cor:image_of_SFT_by_Phi} and Proposition \ref{prop:higher_block_shift_SFT-2}.(C2). Although Proposition \ref{prop:higher_block_shift_SFT-2}.(C2) only had an additional condition added, and its proof remained the same, the loss of generality led to substantial changes in the statement and generality of Corollary \ref{cor:higher_block_shift_Markov_shift}, which was used in most of the theorems in Section 7. In order to preserve the results of Section 7, we added Corollary \ref{cor:higher_block_shift_SFT_on_Z^d-N}, a standard result that was previously contained in the earlier version of Corollary \ref{cor:higher_block_shift_Markov_shift}. It is worth noting that Corollary \ref{cor:higher_block_shift_SFT_on_Z^d-N} is independent of Lemma \ref{lem:SFTs_finite-tersect}, and therefore we can use Theorem \ref{theo:Markov SFT<=>ultragraph presentation} to prove the new version of Lemma \ref{lem:SFTs_finite-tersect}.

We notice that, the change in Proposition \ref{prop:higher_block_shift_SFT-2}.(C2) did not affect the statement or proof of Proposition \ref{prop:sofic and 1-block}. However, Proposition \ref{prop:sofic and 1-block} became more restrictive, leading to significant changes in the statement and generality of Corollary \ref{cor:higher_block_shift_Sofic_on_Z^d-N}. The conditions now imposed in Corollary \ref{cor:higher_block_shift_Sofic_on_Z^d-N} result in a lack of generality in Theorem \ref{theo:wsofic=>graph shift conjugacy} and Theorem \ref{theo:sofic-follower_set_graph-1}. On the other hand, although Corollary \ref{cor:higher_block_shift_Sofic_on_Z^d-N} was also used in theorems \ref{theo:wsofic<=>graph presentation} and \ref{theo:wsofic=>graph presentation-2}, and  (indirectly) in Theorem \ref{theo:sofic-follower_set_graph-2}, its use in those results were restricted to a particular case that now is included in the new Corollary \ref{cor:higher_block_shift_SFT_on_Z^d-N}.

Finally, observe that Theorem \ref{theo:wsofic=>graph shift conjugacy}, which was previously used to prove Claim \ref{claim:FDS-not_weakly_sofic_shifts}, can no longer be used for this purpose.
To preserve Claim \ref{claim:FDS-not_weakly_sofic_shifts}, we introduced two other new independent results: Proposition \ref{prop:higher_block_shift_Sofic_on_Z^d-N} and Theorem \ref{Theo:general_graphs_for_wsofics}.\\

All the changes made with respect to the original article will be presented in the following sections and are highlighted with a gray box. The diagram illustrating the relationships between the results after the corrections introduced in this article is presented in Figure 2.E.

\begin{equation*}\hspace{-2cm}
\xymatrix{
*+[F-:<10pt>]{Corollary\ 5.E }\ar[rd]\ar[r]   &    *+[F-:<10pt>]{Proposition\ 6.E }\ar[r]          &    *+[F-:<10pt>]{Theorem\ 7.E }\ar[r]     &     *+[F]{ Claim\ 9.4 } \\
    &     { \begin{array}{c}Section\ 7\\
   (except\ Theorem\ 7.8)\end{array}}\ar[d]^{Theo.\ 7.3}         &   &    \\
    &     *+[F]{ Lemma\ 5.3 }\ar[r]\ar[ld]      &   *+[F]{Corollary\ 5.4}      &      \\
 *+ [F]{Corollary\ 5.10}  &  *+[F]{Proposition\ 5.8}\ar[r] \ar[l]_{item\ (C2)}       &  *+[F]{Proposition\ 6.5 }\ar[d] &    \\
     &          &  *+[F]{ Corollary\ 6.6 }\ar[d]\ar[rd]   &   \\
     &               &  *+[F]{ Theorem\ 7.4 } \ar[r]        &   *+[F]{ Theorem\ 7.9}
  }
  \end{equation*}
  \captionsetup{labelformat=empty} 

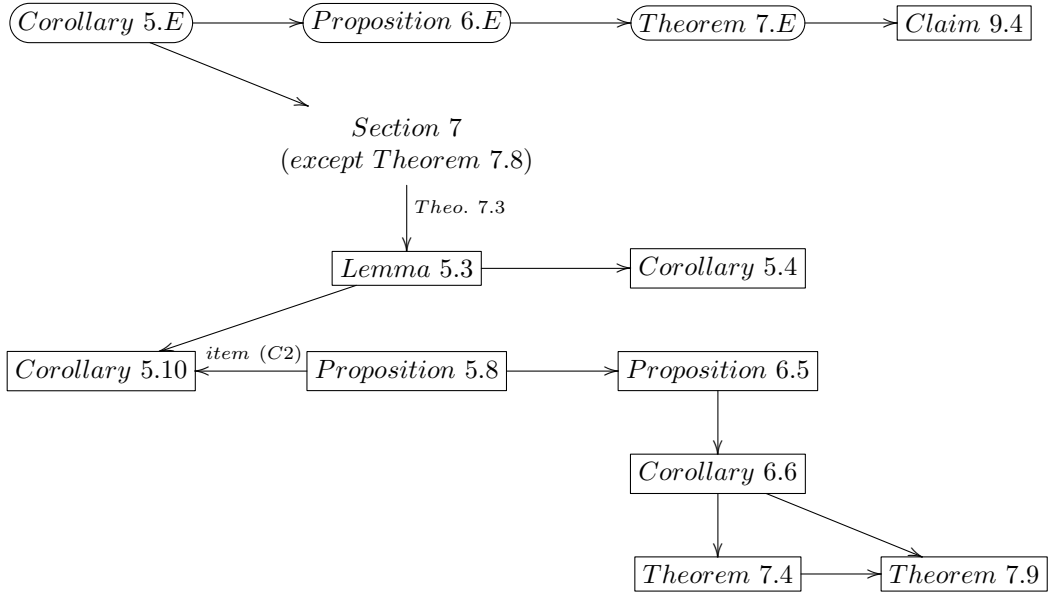
\captionof{figure}{\textbf{Figure 2.E:} The scheme above presents the relationships between the results in the present article. A box around a result indicates that the result was affected by the change in Lemma 5.3. Most of these results have changes in their statements and/or proofs and have become less general (with the exceptions of Proposition 6.5, which retains its statement and proof but is now less general, and Claim 9.4, which has a new proof but retains its generality). This corrected version of the article includes three new results: Corollary \ref{cor:higher_block_shift_SFT_on_Z^d-N}, Proposition \ref{prop:higher_block_shift_Sofic_on_Z^d-N}, and Theorem \ref{Theo:general_graphs_for_wsofics}. These new results are used to provide proofs for most of the theorems in Section 7, as well as for Claim 9.4, which is now independent of Theorem 7.4.}

\vspace{0.75cm}

I am very grateful to Maira Aranguren (UNEXPO, Venezuela), Jorge Campos (UNEXPO, Venezuela), and Neptalí Romero (UCLA, Venezuela) for their careful reading of the original article and for kindly pointing out the flaw in the original version of the proof of Lemma \ref{lem:SFTs_finite-tersect}.

\setcounter{section}{4}
\section*{Changes in Section \ref{sec:SFTs}: Shifts of finite type}

\setcounter{section}{5}
\setcounter{theo}{2}

Below, we present the correct statement and proof of Lemma \ref{lem:SFTs_finite-tersect}, which is considerably more restrictive than the original version.

\begin{lem}\label{lem:SFTs_finite-tersect} If $\NZ$ is \!\!\colorbox{gray!7}{$\N$ or $\Z$ with the usual addition}\!\!, then any \!\!\colorbox{gray!7}{Markov shift (that is, 1-step SFT)}\!\! $\Lambda\subset A^\NZ$  has the property of nested cylinders.
\end{lem}

\begin{tcolorbox}[highlightstyle]
\begin{proof}

If $\Lambda \subset A^{\NZ}$ is a Markov shift, then there exists a directed labeled graph $\mathcal{G}=(G,\mathcal{L})$ that presents it (see Section \ref{sec:graphs} for general results and definitions regarding graph presentations). In particular, from Theorem \ref{theo:Markov SFT<=>ultragraph presentation}.\ref{theo:item:Markov T finite}, we can take $\mathcal{G}$ such that for each fixed symbol $a\in A$, there is an unique vertex where all the edges labeled with $a$ reach.

Let $\{W_\ell\}_{\ell \geq 1}$ be a nested sequence of nonempty cylinders of $\Lambda$. Without loss of generality, we may assume that there exists $T := \{g_i\}_{i \geq 1} \subset \NZ$ such that either $g_{\ell+1} < g_i$ for all $i \leq \ell$, or $g_{\ell+1} > g_i$ for all $i \leq \ell$, and that there exists $(a_{g_i})_{i \geq 1} \in A^T$ such that, for all $\ell \geq 1$, we have $W_\ell = [(a_{g_i})_{i \leq \ell}]_{\Lambda}$. For each $a_{g_i}$ we associated the unique vertex $v_{a_{g_i}}$ where edges labeled as $a_{g_i}$ reach. Given $\ell\geq 1$ define $\alpha(\ell):=\min_{i\leq\ell}\{g_i\}$ and $\beta(\ell):=\max_{i\leq\ell}\{g_i\}$. Therefore,  $W_\ell\neq \emptyset$ implies that in $\mathcal{G}$ there exist a finite
 path $\pi^ \ell:=(\pi^\ell_i)_{\alpha(\ell)+1\leq i\leq\beta(\ell)}$ which:
 \begin{enumerate}
 \item starts at the vertex $v_{\alpha(\ell)}$ and ends at the vertex $v_{\beta(\ell)}$;
 \item for all $\alpha(\ell)<g_i\leq \beta(\ell)$, $\pi^\ell_{g_i}$ is an edge labeled with $a_{g_i}$ (and thus $\pi^\ell_{g_i}$ is and edge reaching the vertex  $v_{a_{g_i}}$).
 \end{enumerate}

Hence, since for all $\ell\geq 1$ we have either $g_{\ell+1}< \alpha(\ell)$ or  $g_{\ell+1}> \beta(\ell)$, we can take $\pi^{\ell+1}$ a path which contains the path $\pi^{\ell}$. Therefore, the sequence $(\pi^\ell)_{\ell\geq 1}$, of paths in $\mathcal{G}$, grows to a (maybe one-sided) infinite path $\pi=(\pi_i)_{i\geq 1}$ in $\mathcal{G}$. Now, observe that, if $\inf_{\ell\geq 1}\{\alpha(\ell)\}=-\infty$, then $\mathcal{L}(\pi_{g_i})=a_{g_i}$ for all $i\geq 1$, and there exists $\x\in \Lambda$, presented (maybe partially) by the path $\pi$, is such that  $x_{g_i}=a_{g_i}$ for all $i\geq 1$, that is, $\x\in\bigcap_{\ell\geq 1}W_\ell$. On the other hand, if $\inf_{\ell\geq 1}\{\alpha(\ell)\}=g_k$ for some $k\geq 1$, then $\pi$ starts at the vertex $v_{a_{g_k}}$. Therefore, we can define the path $\tilde{\pi}=(\tilde\pi_i)_{i\geq 0}$, where $\mathcal{L}(\tilde\pi_0):=a_{g_k}$ and $\tilde\pi_i:=\pi_i$  for all $i\geq 1$. Therefore, $\mathcal{L}(\tilde\pi_{g_i})=a_{g_i}$ for all $i\geq 0$, and so there exists $\x\in\Lambda$ such that $\x\in\bigcap_{\ell\geq 1}W_\ell$.

\end{proof}

\end{tcolorbox}

Therefore, as direct consequence of  Theorem  \ref{theo: image_of_shifts-intersection} and Lemma \ref{lem:SFTs_finite-tersect}, we have:

\begin{cor}\label{cor:image_of_SFT_by_Phi} Suppose $\NZ$ is \!\!\colorbox{gray!7}{$\N$ or $\Z$ with the usual addition}\!\!, and let $\Lambda\subset A^\NZ$ be a \!\!\colorbox{gray!7}{Markov shift}\!\! and $\Phi:\Lambda\to B^\NZ$ be a locally finite-to-one generalized sliding block. Then $\Phi(\Lambda)$ is a shift space.
\end{cor}

\qed\\

Since Lemma \ref{lem:SFTs_finite-tersect} is now less general than its earlier version, item (C2) in Proposition \ref{prop:higher_block_shift_SFT-2} requires an additional condition (its proof remains unchanged).

\setcounter{section}{5}
\setcounter{theo}{7}
\begin{prop}\label{prop:higher_block_shift_SFT-2} Let $\Lambda\subset A^\NZ$ be an SFT  and  $ \bfN$ be a partition of $\Lambda$ by cylinders with $\bfM_\bfN$ finite. If at least one of the following conditions holds:
\begin{enumerate}

\item[(C1)] $\Lambda=A^\NZ$ and $\NZ$ can be extended to a group $\mathbb{G}$ with the property that for all $g,h\in\NZ$ we have $g^{-1}h\in\NZ$ or $h^{-1}g\in\NZ$;

\item[(C2)] $\NZ$ can be extended to a countable group $\mathbb{G}$ with the property that for all $g,h\in\NZ$ we have $g^{-1}h\in\NZ$ or $h^{-1}g\in\NZ$\colorbox{gray!7}{\!, and $\Lambda$ has the property of nested cylinders}\!\!;

\item[(C3)] $\NZ$ is a group and $\bigcap_{M\in\bfM_\bfN}M\neq\emptyset$;

\item[(C4)] $1\in\bigcap_{M\in\bfM_\bfN}M$;

\end{enumerate}
then  $\Lambda^{[\bfN]}$ is also an SFT.

Conversely, if condition (C3) or condition (C4) holds, and $\Lambda^{[\bfN]}$ is an SFT, then $\Lambda$ is also an SFT.

\end{prop}

\qed\\

The original statement of Corollary \ref{cor:higher_block_shift_Markov_shift} was very general and included several cases. Here, we have split this corollary into two parts: a corrected version of Corollary \ref{cor:higher_block_shift_Markov_shift} (which, like its earlier version, is a direct consequence of Proposition \ref{prop:higher_block_shift_SFT-2}); and the new Corollary \ref{cor:higher_block_shift_SFT_on_Z^d-N} (which is independent of Lemma \ref{lem:SFTs_finite-tersect} and its consequences).

\begin{tcolorbox}[highlightstyle]
\setcounter{theo}{9}
\begin{cor}\label{cor:higher_block_shift_Markov_shift}
If $\Lambda\subset A^\NZ$ is a Markov shift, $\NZ$ is $\N$ or $\Z$ with the usual addition, and  $ \bfN$ is a partition of $\Lambda$ by cylinders with $\bfM_\bfN$ finite, then $\Lambda^{[\bfN]}$ is an SFT.
\end{cor}
\end{tcolorbox}

\qed\\

Note that $\Lambda^{[\bfN]}$ in Corollary \ref{cor:higher_block_shift_Markov_shift} is not necessarily a Markov shift. For example, if $\Lambda$ is a Markov shift and $\bfN$ is the family of all cylinders defined on the coordinate 0 and some fixed coordinate $m$, then 
$\Lambda^{[\bfN]}$ will be an $m$-step shift.\\

For classical higher block codes in shift spaces on the classical lattices we have the following result:

\begin{tcolorbox}[highlightstyle]
\renewcommand\thetheo{5.E}
\begin{cor}\label{cor:higher_block_shift_SFT_on_Z^d-N} Let $\NZ$ be the lattice $\N$ or $\Z$ with the usual addition, and for $N\geq 1$ let $\Phi^{[N]}:\Lambda\to \Lambda^{[N]}$ be the $N^{th}$-higher block code given by $\Phi^{[N]}\big((x_i)_{i\in\NZ}\big):=(x_i...x_{i+N-1})_{i\in\NZ}$. It follows that $\Lambda\subset A^\NZ$ is an SFT if and only if  $\Lambda^{[N]}$ is an SFT. In particular, if $\Lambda$ is a $m$-step shift and  $N\geq m$, then $\Lambda^{[N]}$ is a Markov shift.
\end{cor}

\end{tcolorbox}

\begin{tcolorbox}[highlightstyle]
\begin{proof} Just observe that the inverse of 
$\Phi^{[N]}$  is given by
$\left(\Phi^{[N]}\right)^{-1}\big((x_i...x_{i+N-1})_{i\in\NZ}\big)=(x_i)_{i\in\NZ}$, and thus, as $\Phi^{[N]}$ as its inverse are uniform SBCs. Therefore, the proof follows from Theorem \ref{theo:uniform_conjugacy-SFT} (note that this result could be directed proved by using the same arguments that were used in \cite[Theorem 2.1.10.]{LindMarcus}).\\

If $\Lambda$ is an $m$-step shift and $N\geq m$, then a sequence $\y\in (A^\NZ)^{[N]}$ will belong to $\Lambda^{[N]}$ if and only if for each $y_iy_{i+1}= (a_ia_{i+1}...a_{i+N-1})\ (a_{i+1}...a_{i+N-1}a_{i+N})$ we have $a_ia_{i+1}...a_{i+N-1}a_{i+N}\in W_N(\Lambda)$, which means $\Lambda^{[N]}$ is a Markov shift.

\end{proof}

\end{tcolorbox}

\section*{Changes in Section \ref{sec:sofic}: Sofic shifts}

\setcounter{section}{6}
\setcounter{theo}{5}

\begin{tcolorbox}[highlightstyle]
\begin{cor}\label{cor:higher_block_shift_Sofic_on_Z^d-N} Let $\NZ$ be the lattice $\N$ or the lattice $\Z$ with the usual addition.   If $\Omega\subset A^\NZ$ is the image of an shift of finite type $\Gamma$ through a locally finite-to-one sliding block code $\Psi:\Gamma\to\Omega$ with local rule $\psi:W_{\{1\}}(\Gamma)\to A$, then $\Omega$ is weakly sofic shift.  In particular, if $\Psi$ is locally bounded finite-to-one of order $k$, then $\Omega$ is a sofic shift of order $k$.

Conversely, suppose that $\Omega\subset A^\NZ$ is a weakly sofic shift, where $\Lambda\subset B^\NZ$ is an SFT and $\Phi :\Lambda\to\Omega$ is an onto locally finite-to-one SBC, and let $\bfN$ be the partition of $\Lambda$ by the cylinders used in the finitely defined sets that define $\Phi$. If any of the following conditions holds:

\begin{enumerate}

\item[(D1)] $\Lambda$ is a Markov shift

\item[(D2)] $\Lambda$ has the property of nested cylinders

\item[(D3)] $1\in\bigcap_{M\in\bfM_\bfN}M$ (if  $\NZ=\Z$, then it is sufficient that $\bigcap_{M\in\bfM_\bfN}M\neq\emptyset$);

\end{enumerate}
 then there exist a shift of finite type $\Gamma\subset E^\NZ$ and a locally finite-to-one sliding block code $\Psi:\Gamma\to\Omega$ with local rule $\psi:W_{\{1\}}(\Gamma)\to A$ such that $\Omega=\Psi(\Gamma)$. In particular, if $\Omega$ is a sofic shift of order $k$, then $\Psi$ is locally bounded finite-to-one of order $k$.

\end{cor}

\end{tcolorbox}

\begin{proof}
Suppose $\Omega=\Psi(\Gamma)$ for some shift of finite type $\Gamma$ and for some locally finite-to-one sliding block code $\Psi:\Gamma\to\Omega$ with local rule $\psi:W_{\{1\}}(\Gamma)\to A$. Hence, condition (A2) of Theorem \ref{theo: image_of_shifts} holds, and then $\Omega$ is a shift space. Therefore, by definition, it is a weakly sofic shift. In particular, if $\Psi$ is locally bounded finite-to-one of order $k$, then $\Omega$ is a sofic shift of order $k$.\\

\begin{tcolorbox}[highlightstyle]
Conversely, if $\Omega=\Phi(\Lambda)$ for some SFT $\Lambda$ and a locally finite-to-one SBC $\Phi :\Lambda\to\Omega$ such that any of (D1)-(D3) holds, then, since $\NZ=\N$ or $\Z$, and $\bfM_\bfN$ is finite, it follows that some of conditions (C1)-(C4) holds. Therefore, Proposition \ref{prop:sofic and 1-block} ensures that there exist a shift of finite type $\Gamma\subset E^\NZ$ and a locally finite-to-one sliding block code $\Psi:\Gamma\to\Omega$ with local rule $\psi:W_{\{1\}}(\Gamma)\to A$ such that $\Omega=\Psi(\Gamma)$.
In particular, if $\Omega$ is a sofic shift of order $k$, then  $\Psi$ is locally bounded finite-to-one of order $k$.
\end{tcolorbox}

\end{proof}

The next proposition provides a generalization of Corollary \ref{cor:higher_block_shift_Sofic_on_Z^d-N}. However, in the absence of conditions (D1)–(D3), we cannot guarantee that the resulting 1-block code will be locally finite-to-one. This proposition will be used to prove Theorem \ref{Theo:general_graphs_for_wsofics}, which, in turn, is used in the new proof of Claim \ref{claim:FDS-not_weakly_sofic_shifts}.

\begin{tcolorbox}[highlightstyle]
\renewcommand\thetheo{6.E}
\begin{prop}\label{prop:higher_block_shift_Sofic_on_Z^d-N} Let $\NZ$ be the lattice $\N$ or the lattice $\Z$ with the usual addition.   If $\Omega\subset A^\NZ$ is the image of a shift of finite type through a locally finite-to-one sliding block code, then there exist a Markov shift $\Gamma\subset E^\NZ$ and a sliding block code $\Psi:\Gamma\to\Omega$ with local rule $\psi:W_{1}(\Gamma)\to A$ such that $\Omega=\Psi(\Gamma)$. In particular, $W_1(\Gamma)$ is finite if and only if $W_1(\Omega)$ is finite, and $W_1(\Gamma)$ has the same cardinality of $W_1(\Omega)$ whenever $W_1(\Omega)$ is infinite.
\end{prop}
\end{tcolorbox}

\begin{tcolorbox}[highlightstyle]
\begin{proof}

Let $\Omega=\Phi(\Lambda)$ for some SFT $\Lambda\subset B^\NZ$ and some  locally finite-to-one sliding block code $\Phi :\Lambda\to\Omega$. Suppose $\Lambda$ is an $m$-step shift. Let $\ell,r\geq 0$ be the smallest numbers such that we can write the local rule of $\Phi$ as $\phi:W_{\bar m}(\Lambda)\to A$ with $\bar m:=|r-\ell|+1$ and $\Phi(\x)=\Big(\phi(x_{i-\ell}...x_i...x_{i+r})\Big)_{i\in\NZ}$ for all $\x\in\Lambda$ (which exists since $\Phi$ is an SBC).

Take $n$ such that $2n\geq\max\{m,\bar m\}$, and set $N:=2n+1$. Consider the $N^{th}$-higher block code $\Phi^{[N]}:\Lambda\to\Lambda^{[N]}$ given by $\Phi^{[N]}\big((x_i)_{i\in\NZ}\big)=(x_{i-n}...x_i...x_{i+n})_{i\in\NZ}$. By adapting the proof of Corollary \ref{cor:higher_block_shift_SFT_on_Z^d-N} we get that $\Gamma:=\Lambda^{[N]}$ is a Markov shift. Finally, we define $\Psi:\Gamma\to\Omega$ with local rule  $\psi:W_{1}(\Gamma)\to A$ given by $$\psi(a_{-n}...a_0...a_{n}):=\phi(a_{i-\ell}...a_0...a_{i+r}).$$

Observe that since $\Phi$ is assumed to be locally finite-to-one, $W_1(\Lambda)$ and $W_1(\Gamma)$ are finite if and only if $W_1(\Omega)$ is finite, and they have the same cardinality as $W_1(\Omega)$ when the latter is infinite.

\end{proof}

\end{tcolorbox}

We remark that, in the above proposition, it is not required that $\Omega$ is a shift space (that is, it could be a non-closed set). Furthermore, even if $\Omega$ is a weakly sofic shift, that is, even if $\Omega$ is a shift space and $\Omega=\Phi(\Lambda)$ for some locally finite-to-one SBC $\Phi$, the SBC $\Psi$ found might be not locally finite-to-one (for example, if we take $\Omega=\Lambda:=A^\Z$, $\Phi:\Lambda\to\Omega$ as $\Phi=\s$, and $\Gamma=\Lambda^{[2]}$, the corespondent 1-block map $\Psi$ is not locally finite to one). We also remark that the fact that the cardinality of $W_1(\Gamma)$ is related to that of $W_1(\Omega)$ is the main result of Proposition \ref{prop:higher_block_shift_Sofic_on_Z^d-N}, since any shift space on the lattice $\N$ or $\Z$ is always the image of a Markov shift $\Gamma \subset \Omega^\NZ$ under a sliding block code that is not locally finite-to-one (see the discussion preceding Theorem \ref{Theo:general_graphs_for_wsofics}).

\section*{Changes in Section \ref{sec:graphs}: Graph presentations of weakly sofic shifts on the lattices $\N$ and $\Z$}

We remark that, in the original article, most of the results in Section \ref{sec:graphs} relied on Corollary \ref{cor:higher_block_shift_Markov_shift} to ensure that the $N^{th}$-higher block code of an SFT is also an SFT. Now, this result no longer follows from Corollary \ref{cor:higher_block_shift_Markov_shift}, but rather from Corollary \ref{cor:higher_block_shift_SFT_on_Z^d-N}. Thus, in the paragraph before Theorem \ref{theo:SFT<->edge shift conjugacy}, the sentence\\

{\em
``In particular, if $\Lambda$ is an SFT, then  $\Lambda^{[N]}$ is also an SFT (Corollary \ref{cor:higher_block_shift_Markov_shift} or it can also be proved using the same arguments used in \cite[Theorem 2.1.10.]{LindMarcus}).’’\\
}

shall be replaced by\\

{\em 
``In particular, \!\!\colorbox{gray!7}{$\Lambda$ is an SFT if and only if $\Lambda^{[N]}$ is an SFT}\!\! (\colorbox{gray!7}{\!\!Corollary \ref{cor:higher_block_shift_SFT_on_Z^d-N}}\!\! or it can also be proved using the same arguments used in \cite[Theorem 2.1.10.]{LindMarcus}).’’\\
}

We notice that Theorem \ref{theo:SFT<->edge shift conjugacy} and Theorem \ref{theo:Markov SFT<=>ultragraph presentation} ensure that shifts of finite type can always be presented by a directed labeled graph, regardless of whether the lattice is $\N$ or $\Z$. Furthermore, Example \ref{ex:free_context} shows that some shift spaces on the lattice $\Z$, which are not of finite type, can also be presented by directed labeled graphs. Although Theorem \ref{theo:motivation} can only be applied to find labeled graphs presenting shift spaces on the lattice $\N$, for shift spaces on the lattice $\Z$ one could always define a labeled graph with as many disjoint paths as there are orbits in the shift. However, such a graph would have an uncountable number of vertices and edges, even for most shift spaces over a finite alphabet.  Therefore, any interesting extension of Theorem \ref{theo:motivation} to the lattice $\Z$ should encapsulate only graphs whose vertex and edge set cardinalities are related to the cardinality of the alphabet of $\Lambda$. In fact, there exist shift spaces on the lattice $\Z$ and with countable alphabets, for which is not possible to find a graph presentation with a countable number of vertices and edges (see Claim \ref{claim:FDS-not_weakly_sofic_shifts} for an example).
The following theorem gives an extension of Theorem \ref{theo:motivation} for some class of $\s$-invariant subsets.

\begin{tcolorbox}[highlightstyle]
\renewcommand\thetheo{7.E}
\begin{theo}\label{Theo:general_graphs_for_wsofics}  Let $\NZ$ be the lattice $\N$ or $\Z$  with the usual addition, and suppose $\Omega\subset A^\NZ$ is the image of a shift of finite type through a locally finite-to-one sliding block code. If $W_1(\Omega)$ is finite, then $\Omega$ can be presented by a labeled graph whose set of edges and set of vertices are finite. If $W_1(\Omega)$ is infinite, then $\Omega$ can be presented by a labeled graph whose set of edges has the same cardinality of  $W_1(\Omega)$, and whose set of vertices has cardinality no greater than the cardinality of  $W_1(\Omega)$.
\end{theo}
\end{tcolorbox}

\begin{tcolorbox}[highlightstyle]
\begin{proof}

If $\Omega\subset A^\NZ$ is the image of a shift of finite type through a locally finite-to-one sliding block code, then, from Proposition \ref{prop:higher_block_shift_Sofic_on_Z^d-N}, there exist a Markov shift $\Gamma\subset E^\NZ$ and a sliding block code $\Psi:\Gamma\to\Omega$ with local rule $\psi:W_{1}(\Gamma)\to A$ such that $\Omega=\Psi(\Gamma)$. In particular, $W_1(\Gamma)$ is finite whenever $W_1(\Omega)$ is finite, and it has the same cardinality as $W_1(\Omega)$ when the latter is infinite. Since, $\Gamma$ is a Markov shift, from Theorem \ref{theo:Markov SFT<=>ultragraph presentation}, there exists a labeled graph $\mathcal{G}=(G,\mathcal{L})$ that presents $\Gamma$ whose set of edges is $\mathcal{E}_G=W_2(\Gamma)$ and set of vertices is $\mathcal{V}_G=W_1(\Gamma)$. Hence, the labeled graph $\widehat{\mathcal{G}}:=(G,\psi\circ\mathcal{L})$ is a presentation of $\Omega$ that satisfies the statement of the theorem.\\

Note that the cardinality of $W_1(\Omega)$ is both, a lower bound for the cardinality of the edge set of any presentation of $\Omega$ (since each symbol in $W_1(\Omega)$ should be associated at least to one edge in any graph presentation of $\Omega$), and an upper bound for the minimum quantity of vertices for which we can find a presentation of $\Omega$ (since the previous construction ensures that exists one presentation using $W_1(\Gamma)$ as vertex set).

\end{proof}
\end{tcolorbox}

Observe that, it was not required in the above theorem that $\Lambda$ is a shift space. Furthermore, even if $\Omega$ was a weakly sofic shift, the labeled graph given by the theorem may contain labels that appear infinitely many times, possibly as outgoing edges from, or ingoing edges to, infinitely many vertices (because $\psi$ might be not finite-to-one). Theorems \ref{theo:wsofic=>graph shift conjugacy}-\ref{theo:sofic-follower_set_graph-2} in \cite{Sobottka2022} provided conditions under which we can establish several connections between weakly sofic shift spaces and directed labeled graphs with special features.\\

The following two theorems now include an additional condition, which makes them less general. Their proofs follow as in the original article.

\setcounter{section}{7}
\setcounter{theo}{3}

\begin{theo}\label{theo:wsofic=>graph shift conjugacy} If $\Lambda\subset A^\NZ$ is a weakly sofic shift \!\!\colorbox{gray!7}{as in Corollary \ref{cor:higher_block_shift_Sofic_on_Z^d-N}}\!\!, then there exists $M\in\N$ such that
\begin{enumerate}

\item\label{theo:wsofic=>graph shift conjugacy_1} $\Lambda^{[M+1]}=X_\mathcal{\hat G}$ for some directed labeled graph $\mathcal{\hat G}$ where each label is used just finitely many times. In the particular case of $\Lambda$ being a sofic shift of order $k$, then each label is used at most $k^{M+1}$ times in $\mathcal{\hat G}$.

\item\label{theo:wsofic=>graph shift conjugacy_2} $\Lambda=X_\mathcal{G}=X_\mathcal{H}$ where $\mathcal{G}$ and $\mathcal{H}$ are directed labeled graphs such that for all $w\in W_m(\Lambda)$ with $m\geq M$ we have $|I_\mathcal{G}(w)|<\infty$ and $|T_\mathcal{H}(w)|<\infty$. In the particular case of $\Lambda$ being a sofic shift of order $k$, then for all $w\in W_m(\Lambda)$ with $m\geq M$ we have $|I_\mathcal{G}(w)|\leq k^m$ and $|T_\mathcal{H}(w)|\leq k^m$.

\end{enumerate}

\end{theo}

\qed\\

\setcounter{section}{7}
\setcounter{theo}{8}

\begin{theo}\label{theo:sofic-follower_set_graph-1} Let $\Lambda\subset A^\NZ$ be a shift space and let $\mathcal{G}_\mathcal{F}$ be its follower set graph. If $\Lambda$ is a weakly sofic shift \!\!\colorbox{gray!7}{as in Corollary \ref{cor:higher_block_shift_Sofic_on_Z^d-N}}\!\!, then there exists $M\geq 0$ such that for all $m\geq M$ and $w\in W_m(\Lambda)$ we have $|T_{\mathcal{G}_\mathcal{F}}(w)|<\infty$. In particular, if $\Lambda$ is a sofic shift of order $k$, then for
 all $m\geq M$ and $w\in W_m(\Lambda)$ we have $|T_{\mathcal{G}_\mathcal{F}}(w)|\leq 2^{k^m}-1$. Moreover, if $\Lambda$ is an $M$-step SFT, then  for
 all $m\geq M$ and $w\in W_m(\Lambda)$ we have $|T_{\mathcal{G}_\mathcal{F}}(w)|= 1$.
\end{theo}

\qed\\

\setcounter{section}{8}

\section*{Changes in Section \ref{sec:relationships}: Relationship between shift spaces}

The original proof of Claim \ref{claim:FDS-not_weakly_sofic_shifts}, given in \cite{Sobottka2022}, used Theorem \ref{theo:wsofic=>graph shift conjugacy}.\ref{theo:wsofic=>graph shift conjugacy_1} to show that the variable length shift of Example \ref{ex:SVL1} is not a weakly sofic shift. However, since Theorem \ref{theo:wsofic=>graph shift conjugacy} is now valid only for weakly sofic shifts satisfying some of the conditions (D1)-(D3) of Corollary \ref{cor:higher_block_shift_Sofic_on_Z^d-N}, the original argument used in the proof of Claim \ref{claim:FDS-not_weakly_sofic_shifts} just proved that the shift of the Example \ref{ex:SVL1} is not a weakly sofic shift satisfying (D1)-(D3).\\

To give a complete proof for Claim \ref{claim:FDS-not_weakly_sofic_shifts}, we will use other SVL, and Theorem \ref{Theo:general_graphs_for_wsofics}.

\setcounter{section}{9}
\setcounter{theo}{3}

\begin{claim}\label{claim:FDS-not_weakly_sofic_shifts} There exist FDSs that are not weakly sofic shifts.
\end{claim}

\begin{tcolorbox}[highlightstyle]

Clearly an FDS which is not a weakly sofic shift shall be an SVL. Let $A=\N$ and consider the full shift $A^\Z$. Let $\bfN$ be the partition of $A^\Z$ which contains the cylinder in the form $[x_0=0]$, and $[x_0=n,x_n=k]$ for all $n,k\in A$, and take the SVL $(A^\Z)^{[\bfN]}$. Note that $\y\in(A^\Z)^{[\bfN]}$ can be written as

$$\y=\left(\begin{array}{l}x_i\\x_{i+x_i}\end{array}\right)_{i\in\Z},$$
for some $\x\in A^\Z$.

Suppose, by contradiction, that $(A^\Z)^{[\bfN]}$ is a weakly sofic shift. Let $\mathcal{G}$ be any directed labeled graph that presents $(A^\Z)^{[\bfN]}$ (whose existence is guaranteed by Theorem \ref{Theo:general_graphs_for_wsofics}). Given $\mathbf{a}:=(a_i)_{i\leq-1}\in A^{-\N*}$, consider the left infinity sequence 
$$\y_{\mathbf{a}}^-:=\left(\begin{array}{l}y_i\\y_{i+y_i}\end{array}\right)_{i\leq-1}=\left(\begin{array}{l}2|i|-1\\a_i\end{array}\right)_{i\leq-1}.$$

Note that, $\y_{\mathbf{a}}^{-}$ can be uniquely extended to a sequence $\y_{\mathbf{a}}\in(A^\Z)^{[\bfN]}$. This means that any left sided infinite path in $\mathcal{G}$ that represents $\y_{\mathbf{a}}^{-}$ will end in a same vertex $v_{\mathbf{a}}$ from which there is only one possible way to follow. Furthermore, given $\mathbf{a},\mathbf{b}\in A^{-\N*}$, $\mathbf{a}\neq\mathbf{b}$, the continuation of $\y_{\mathbf{a}}^{-}$ is different of the continuation of $\y_{\mathbf{b}}^{-}$, and then $v_{\mathbf{a}}\neq v_{\mathbf{b}}$. Hence, since there are uncountable many choices of $\mathbf{a}$, it means that $\mathcal{G}$ has uncountable many vertices $v_{\mathbf{a}}$, contradicting Theorem \ref{Theo:general_graphs_for_wsofics} that should exist a labeled graph whose cardinality of the vertex set is not greater than the cardinality of $W_1((A^\Z)^{[\bfN]})$.

\end{tcolorbox}

Hence, Figure 7 in \cite{Sobottka2022}, shall be replaced by the following one:

\begin{figure}[H]
\centering
\includegraphics[width=0.8\linewidth=1.0]{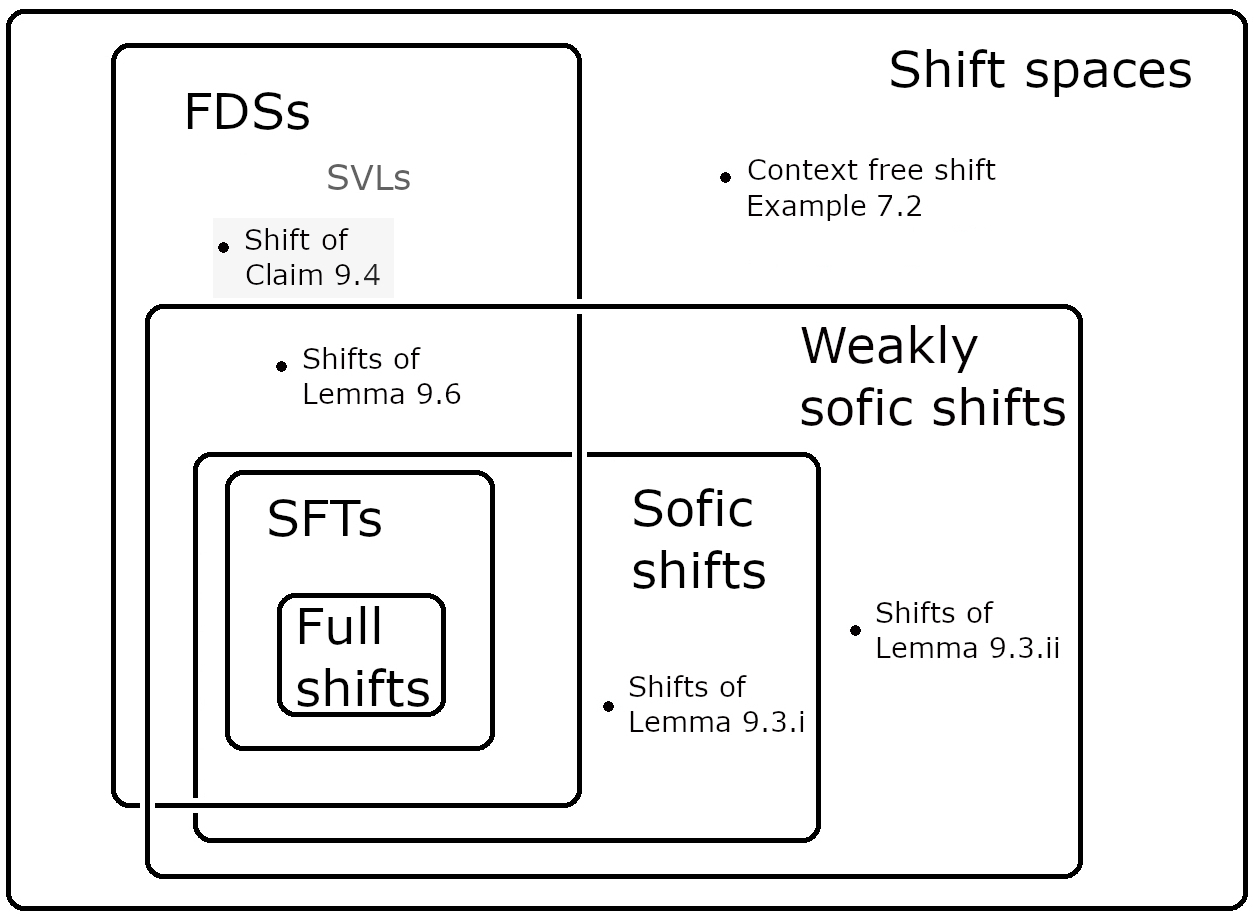}
\captionof{figure}{\textbf{Figure 7:} Relationship between classes of shift spaces.}\label{Fig_shift_spaces_relationship}
\end{figure}

\section*{Acknowledgments}

The author acknowledges the assistance of ChatGPT (OpenAI) in improving the grammatical clarity of the manuscript.


\end{document}